\documentclass[11pt]{article}
\usepackage{amssymb,amsfonts,amsthm,amsmath}
\usepackage{color}
\usepackage{footmisc}

\newtheorem{thm}{Theorem}[section]

\newtheorem{lem}[thm]{Lemma}
\newtheorem{prop}[thm]{Proposition}

\newtheorem{defn}[thm]{Definition}
\newtheorem{rem}[thm]{Remark}

\def\eop{\hfill\rule{2.5mm}{2.5mm}}
\def\pf{\par\smallbreak\noindent {\bf Proof.} \ }

\begin{document}

\title{Reproducing Kernel Hilbert Spaces and entropy Kolmogorov numbers on compact Lie Groups
}

\author{Z. Avetisyan\footnotemark[1]
\and
K. Gonzalez\footnotemark[2]\;\,\footnotemark[3]
\and M. Ruzhansky\footnotemark[1]\;\,\footnotemark[4]}

\renewcommand{\thefootnote}{\fnsymbol{footnote}}

\footnotetext[1]{Ghent University, Belgium. Department of Mathematics: Analysis, Logic and Discrete Mathematics}
\footnotetext[2]{Yerevan State University, Armenia. Faculty of Mathematics and Mechanics} \footnotetext[3]{Weizmann Institute of Science, Israel. Faculty of Mathematics and Computer Science}
\footnotetext[4]{Queen Mary University London, UK. Department of Mathematics}

\maketitle

\begin{abstract}
    On a compact Lie group $G$, we consider the reproducing kernel Hilbert space $\mathcal{H}_K$ associated with the integral kernel $K$ of a left-invariant, positive, symmetric, trace class integral operator on $L^2(G)$. We present lower and upper asymptotic estimates for the entropy Kolmogorov numbers (also called covering numbers) for the embedding of $\mathcal{H}_K$ into the space $C(G)$ of continuous functions on $G$.
\end{abstract}


\thispagestyle{empty}

\section{Introduction}
The Reproducing Kernel Hilbert Space (RKHS) is a tool from functional analysis that sets up a formal way to analyze functions using kernels. An RKHS $\mathcal{H}_K$ is uniquely associated with a symmetric and positive definite kernel function $K(x, y)$, defined on $X\times X$, with $X$ a metric space \cite{aronszajn}. The defining feature of an RKHS is the reproducing property, which states that for any function $f \in \mathcal{H}_K$ and any point $x$ in the domain $X$, $f(x) = \langle f, K(\cdot, x) \rangle_{\mathcal{H}_K}$, linking the function's value to an inner product involving the kernel. This property makes the RKHS a natural setting for machine learning and approximation theory, as it allows for the analysis of functions in a high-dimensional feature space without explicitly computing the feature map, relying instead on the kernel trick \cite{scholkopf_smola_2002, jordao-men}.

While the RKHS naturally arises from kernels that are positive definite, the early stages of kernel machine learning in the 1990s focused heavily on kernels fulfilling the specific conditions outlined in Mercer's theorem \cite{cristianini_shawe-taylor_2000, Mercer1909, Vapnik1998}. In recent years, estimation and learning methods that employ positive definite kernels have become quite popular, since positive definite kernels seem to be the right class of kernels to consider.

In a metric space, the minimal number of balls with a given radius $\epsilon$ required to cover a subset is called the $\epsilon$-covering number of that set. Covering numbers of the unit ball's image belonging to the RKHS via an embedding operator or covering numbers of the embedding operator turned out to be a key tool for evaluating the potential error of statistical algorithms based on observed data \cite{smale,zhou1,zhou}. Developing techniques for obtaining sharp estimates for the covering numbers open the door to resolving many questions in applied mathematical fields. 
Also known as entropy Kolmogorov numbers \cite{kolmogorov}, the covering numbers contribute significantly in areas such as kernel-based learning algorithms and Gaussian process \cite{Quang,ingobook,WSS}.

In 2002, Zhou undertook a rigorous analysis of the covering numbers of the images of ball within RKHS through the embedding into the space of continuous function. By leveraging the relationship between covering numbers and smoothness properties of the Fourier transform of the kernels, Zhou established upper bounds for covering numbers related to Gaussian RKHS and formulated a compelling conjecture regarding their asymptotic behavior when the radius of the balls approaches zero from the right \cite{zhou1}.

Years later, Kuhn in 2011 established the exact asymptotic behavior of these covering numbers. His innovative strategy hinged on an explicit representation of Gaussian RKHS using orthonormal bases \cite{Kuhn}, where the essential tools in his proof were the results of Steinwart et al. \cite{Ingo} on the structure of Gaussian RKHS; in particular, the specific orthonormal bases of these spaces given in \cite{Ingo} were explored.

Inspired by this line, in 2023 sharp estimates for the covering numbers of the embedding $I_W : \mathcal{H}_W \rightarrow C(I)$ of the RKHS $\mathcal{H}_W$ associated with the Weierstrass fractal kernel $W$ into the space of continuous functions $C(I)$ were presented. Here $C(I)$ is the space of continuous real-valued functions on $I\doteq[-1,1]$ endowed with the supremum norm $\|\, \cdot \,\|_{\infty}$  \cite{Gonzalez}. Later, as a generalisation of this work, in 2024 precise estimates for the covering numbers of the embedding of two general classes of RKHS on the unit sphere \cite{Gonzalez2} and on compact two-point homogeneous spaces were provided \cite{gonzalez3}.

The study of the covering numbers in classical RKHS settings has been characterised by the establishment of sharp asymptotic equivalence between their upper and lower bounds. Indeed, when the RKHS comes from a kernel $K$ with $X$ a closed interval or $X$ a compact two-point homogeneous domain of dimension $d$ with the sphere as a particular case with Fourier or Schoenberg coefficients dominated by a geometric progression, the entropy Kolmogorov numbers of the immersion operator $I_K$ has behavior equivalent to $[\ln (1/\epsilon)]^{d+1}$, as $\epsilon\rightarrow 0^+$. This behavior changes dramatically when we consider a kernel $K$ with Fourier coefficients decaying as a harmonic progression \cite{Gonzalez2,gonzalez3}.
 
In this work, we extend the study of the behavior of covering numbers to the setting of compact Lie groups. We start by exploring the matrix form of the unitary representations of compact Lie groups, providing Fourier series expansions for $L^2$-functions. This, combined with the global quantization and symbols of pseudo-differential operators on compact groups (e.g., the Fourier transform of a right convolution kernel) \cite{Michael2}, allows us to explicitly write the kernel $K$ in terms of the symbol, facilitating the study of the properties of $K$. We are mainly interested in studying the RKHS $\mathcal{H}_K$ associated with the kernel $K$, in order to find bounds for the covering numbers of the embedding operator $I_K :\mathcal{H}_K \rightarrow C(G)$, where $C(G)$ is the space of continuous functions on the compact Lie group $G$. Exploring RKHS and covering numbers on compact Lie groups combines tools from abstract harmonic and functional analysis, and is thus a novel subject that promises to open new fronts of study in mathematics.

The paper is organized as follows. In Section 2, we present the main facts about compact Lie groups to be used in the sequel, and we introduce the integral kernel that will be the focus of our study. In Section 3, we characterize the continuity, symmetry, and positivity of the kernel in terms of its symbol. In Section 4, we consider a trace class invariant operator and present the RKHS associated to the kernel presented in the previous section. In Section 5, we define the trace and determinant of finite order of the symbol and present estimates for covering numbers of the embedding operator $I_K$.



\section{Background} In this section, we recall some basic facts of square integrable functions defined over a compact Lie group. Moreover, we introduce the operators from Schatten-von Neumann class and the representation of an integral kernel over a compact Lie group that we will use throughout this work.

Let $L^2\left(G\right)\doteq L^2\left(G, d\mu\right)$ be the Hilbert space of square integrable functions $f: G \rightarrow \mathbb{C}$, with $G$ a compact Lie group of dimension $n$ and $\mu $ the normalized Haar measure of $G$. We consider the norm $\|\cdot\|_2$ induced by the inner product
\begin{equation}\label{inner}
\langle f, g\rangle_2=\int_{G} f(x) \overline{g(x)} d \mu(x), \quad f, g \in L^2\left(G\right).
\end{equation}
Throughout this work we denote $d\mu(x) $ simply by $dx$.

Let $ \hat{G}$ denote the set of equivalence classes of
continuous irreducible unitary representations $\xi$ of $G$. For $[\xi]\in \hat{G}$, by choosing an orthonormal basis in the representation space of $\xi$, we can view $\xi$ as a unitary matrix-valued function $\xi: G \rightarrow \mathbb{C}^{d_{\xi} \times{d_\xi}}$, where $d_\xi$ is the dimension of the
representation space of $\xi$. By the Peter and Weyl Theorem \cite[Theorem 7.5.14]{Michael1} we know that the set 
\[
\mathcal{B} = {\left\{\sqrt{d_\xi} \,\xi_{ij}: (\xi_{ij})_{i,j=1}^{d_\xi} \, , [\xi]\in \hat{G}\right\}},
\] 
is an orthonormal basis for $L^2 (G)$.
Thereby, if $f\in L^2 (G)$, 
\[
f(x)= \sum_{[\xi]\in \hat{G}}\sum_{i,j=1}^{d_\xi} \langle f, \sqrt{d_\xi} \,\xi_{ij}\rangle_2 \, \sqrt{d_\xi}\, \xi_{ij} (x),
\]for almost every $x\in G$.
Arranging the terms and remembering that $\langle f, \xi_{ij}\rangle_2 =(\hat{f} (\xi))_{ji} $
we have 
\begin{align*}
f(x)&= \sum_{[\xi]\in\hat{G}} d_\xi \sum_{j=1}^{d_\xi} \sum_{i=1}^{d_\xi} (\hat{f} (\xi))_{ji} \xi_{ij} (x)\\
    &=\sum_{[\xi]\in\hat{G}} d_\xi \sum_{j=1}^{d_\xi} [\hat{f}(\xi) \,\xi(x)]_{jj}.
\end{align*}
Therefore, as a consequence of the Peter-Weyl theorem, we have the Fourier inversion formula given by
\begin{equation}\label{formula}
f(x)=\sum_{[\xi]\in\hat{G}} d_\xi\, \mathrm{Tr}[\hat{f}(\xi)\, \xi(x)]\,,
\end{equation}
for almost every $x\in G$.
In other words, if $\mathfrak{F}: L^2 (G)\rightarrow \ell^2 (\hat{G}) $ is the group Fourier transform, we also write the inverse Fourier transform as
\begin{equation}\label{formula2}
(\mathfrak{F}^{-1} C)(x)=\sum_{[\xi]\in\hat{G}} d_{\xi} \mathrm{Tr} [C(\xi)\,\xi(x)],
\end{equation}
for almost every $x\in G$ \cite[Theorem 10.3.23]{Michael1}.

We consider a compact operator $T:  L^{2} (G)\rightarrow L^{2} (G)$. Then $T^* T$ is compact, self-adjoint, and non-negative. Hence, we can define the absolute value of $T$ by the equality $|T|=\left(T^* T\right)^{\frac{1}{2}}$. Let $s_k(T)$ be the eigenvalues of $|T|$ with $k=1,2, \ldots$, i.e., the eigenvalues of $ |T|$ with multiplicities counted. The numbers $s_1(T) \geq s_2(T) \geq \cdots \geq s_n(T) \geq \cdots \geq 0$, are called the singular values of $T$. If $0<p<\infty$ and the sequence of singular values is $\ell^p$-summable, then $T$ is said to belong to the \emph{Schatten-von Neumann class} $S_p\left(L^2 (G),L^2 (G)  \right)\doteq S_p (L^2 (G) )$. If $1 \leq p<\infty$, a norm is associated to $S_p\left(L^2 (G) \right)$ by
\[
\|T\|_{S_p (L^2 (G))}\doteq\left(\sum_{k=1}^{\infty}\left(s_k(T)\right)^p\right)^{\frac{1}{p}} .
\]

If $1 \leq p<\infty$ the class $S_p\left(L^2 (G) \right)$ becomes a Banach space endowed with the norm $\|\cdot\|_{S_p (L^2 (G))}$. 
 If $p=\infty$ we define $S_{\infty (L^2 (G))}$ to be the class of compact operators from $L^2 (G) $ into $L^2 (G) $ endowed with the operator norm $\|T\|_{S_{\infty}(L^2 (G))}\doteq\|T\|_{op}$. To further abbreviate the notation, we will write $S_p$ for $S_p (L^2 (G))$. 
 
The Schatten-von Neumann classes are nested, with
\begin{equation}\label{inc}
S_p \subset S_q , \quad\text { if } 1\leq p<q \leq \infty,
\end{equation}
and satisfy the composition embedding 
\begin{equation*}
S_q S_p \subset S_r,
\end{equation*}
where
\begin{equation*}\label{2}
\frac{1}{r}=\frac{1}{p}+\frac{1}{q}, \quad 1\leq p<q \leq \infty .
\end{equation*}

Moreover, for $A\in S_p$ and $B\in S_q$ one has
\begin{equation}\label{1}
\|A B\|_{S_r} \leq\|A\|_{S_p}\|B\|_{S_q},
\end{equation}
and from e.g. \cite[Theorem 3.1]{Delgado}, we have
\begin{equation}\label{Tr}
\|A\|_{S_p}= (\mathrm{Tr}(|A|^p))^{1/p},\quad 1 \leq p < \infty .
\end{equation}
The classes $S_2 $ and $S_1 $ are usually known as the \emph{class of Hilbert-Schmidt} operators and the \emph{trace class}, respectively \cite{Delgado3}. From \cite[Theorem 5]{Ruskai}, if $A$ is a trace class operator then  
\begin{equation}\label{Tr2}
|\mathrm{Tr} A |\leq \|A\|_{S_1}.
\end{equation}

Let us consider $T$ a trace class operator, $D(G)$ the set of the test functions and its dual space $D^{'}(G)$, the space of distributions. Note that since $G$ is compact, we have that $D(G)= C^{\infty} (G)$. 
We remember that the Haar measure is invariant under left and right translations on $G$ and we announce the following statement.

\begin{prop}\cite[Corollary 3.7]{Delgado2}\label{Left-translation}
Let $G$ be a compact Lie group and let the left-invariant operator $T$ be bounded on $L^2 (G)$. Then it is of the form $Tf = f * k$ with $k\in D^{'} (G)$. Moreover, if  $T\in S_2 $, then $k \in L^{2} (G)$.
\end{prop}


In other words, if ${T}$ is invariant under left-translations, then $T:L^2 (G)\rightarrow L^2 (G) $ is a convolution operator with (right convolution) kernel ${k}$ such that  
\begin{equation}\label{T}
Tf(x)=(f*k)(x)=\int_G k(y^{-1} x)f(y)dy,\quad a.e.\,x\in G,
\end{equation}
for all $f\in L^2 (G)$. Since we assumed $T\in S_1$, by the formula (\ref{inc}) we have $T\in S_2$, and from Proposition \ref{Left-translation} we obtain $k\in L^2 (G)$.

Using the formula (\ref{formula}) we can write
\begin{align*}
k(y^{-1} x)&= \sum_{[\xi]\in\hat{G}} d_\xi\, \mathrm{Tr}[\hat{k}(\xi)\, \xi (y^{-1} x)]\\
           &= \sum_{[\xi]\in\hat{G}} d_\xi\, \mathrm{Tr}[ \xi (x) \hat{k}(\xi)\, \xi (y)^* ],
\end{align*}
for almost every $x,y\in G$.
If $K:G\times G \rightarrow \mathbb{C}$ is a kernel defined by $K(x,y)\doteq k(y^{-1} x)$ we have
\[
K(x,y)= \sum_{[\xi]\in\hat{G}} d_\xi\, \mathrm{Tr}[ \xi (x) \hat{k}(\xi)\, \xi (y)^* ],
\]
for almost every $x,y\in G$, and thereby we can write the formula (\ref{T}) as
\begin{equation}\label{T2}
Tf(x)=\int_G K(x,y)f(y)dy,\quad a.e.\, x\in G,
\end{equation}
for all $f\in L^2 (G)$.

We introduce $\sigma_T (\xi)$, the matrix symbol of $T$ of size $d_\xi \times d_\xi$, which we define by $\sigma_T (\xi)\doteq \hat{k}(\xi)$ with $[\xi]\in \hat{G}$. Since $k\in L^2 (G)$, $\sigma_T =\hat{k}\in \ell^2 (\hat{G})$. Therefore, we can write
\begin{equation}\label{Kernel}
K(x,y)= \sum_{[\xi]\in\hat{G}} d_\xi\, \mathrm{Tr}[ \xi (x) \sigma_T (\xi)\, \xi (y)^* ],
\end{equation}
 for almost every $x,y\in G$.


\section{Continuous, symmetric and positive definite kernels on compact Lie groups}

In what follows we are interested in knowing when the kernel $K$ is continuous, symmetric and positive definite. A function $K:G\times G \rightarrow \mathbb{C}$ is a \emph{symmetric kernel} if $K(x,y)=\overline{K(y,x)}$, for almost every $x,y\in G$, and is a \emph{positive definite kernel} if for any $n\in\mathbb{N}$, $x_1,\ldots,x_n \in G$ and $c_1, \ldots,c_n\in \mathbb{C}$ it holds that
\[
\sum_{i=1}^{n}\sum_{j=1}^{n}c_{i}\overline{c_{j}} K(x_{i},x_{j})\geq 0.
\]
The last definition is equivalent to having $z^* (K(x_i , x_j ))_{i,j=1}^{n}\, z \geq 0$, for every nonzero complex column vector $z$ in $\mathbb{C}^n$, where $z^{*}$ denotes the conjugate transpose of $z$.
The integral linear operator $T: L^2 (G)\rightarrow L^2 (G)$ is a \emph{positive definite operator} if $\langle Tf,f\rangle_2 \geq 0$, for all $f\in L^2 (G) $.

\begin{lem}\label{Continuity}
The kernel $K$ defined in (\ref{Kernel}) is continuous whenever $T$ is a trace class operator on $L^2 (G)$.
\end{lem}
\pf
Due to the inequality (\ref{1}) ($r=1$, $p=1$ and $q=\infty$), formule (\ref{Tr2}) and (\ref{Tr}), we obtain
\[
|\mathrm{Tr}{[ \xi(x)\sigma_T (\xi) \xi(y)^* ]}| \leq \| \xi(y^{-1} x) \|_{op} \|\sigma_T (\xi) \|_{S_1} 
=  \mathrm{Tr} (|\sigma_T (\xi)|) ,
\]
for almost every $x,y\in G$, and consequently,
\[
\sum_{[\xi]\in \hat{G}} |d_{\xi} \mathrm{Tr}[\xi(x)\sigma_T (\xi) \xi(y)^* ]| \leq \sum_{[\xi]\in \hat{G}} d_{\xi} \mathrm{Tr} (|\sigma_T (\xi)|).
\]
If  $T$ is a trace class operator, then 
\[
\sum_{[\xi]\in \hat{G}} d_{\xi} \mathrm{Tr} (|\sigma_T (\xi)|)< \infty, 
\]
\cite[Theorem 3.5]{Delgado2}. By the M-test of Weierstrass, the series defining $K(x,y)$ converges uniformly in $G\times G$ almost everywhere, thereby $K(x,y)$ is continuous. 
\eop

\begin{lem}\label{Lm1} The kernel $K$ given by (\ref{Kernel}) is symmetric if and only if the 
matrix symbol $\sigma_T (\xi)$ is a hermitian matrix for all $[\xi]\in\hat{G}$.
\end{lem} 

\pf If $\sigma_T (\xi)^* = \sigma_T (\xi)$ for all $[\xi]\in\hat{G}$, then we have 
\begin{align*}
    \overline{K(y,x)}
=& \sum_{[\xi]\in\hat{G}} d_\xi \,\overline{\mathrm{Tr}[\xi (y) \sigma_T (\xi) \xi (x)^*] }\\
    =& \sum_{[\xi]\in\hat{G}} d_\xi \,\mathrm{Tr}[(\xi (y) \sigma_T (\xi) \xi (x)^*)^*] \\
    =& \sum_{[\xi]\in\hat{G}} d_\xi \,\mathrm{Tr}[\xi (x) \sigma_T (\xi) \xi (y) ^* ]\\
    =& K(x,y),
\end{align*}
for almost every $x,y\in G$. If $K$ is a symmetric kernel we obtain $\sigma_T (\xi)^* = \sigma_T (\xi)$ by the uniqueness of the representation in Fourier series, completing the proof.
\eop

Furthermore, if $T$ is a bounded operator, then $K$ is a symmetric kernel if and only if the operator $T$ is self-adjoint \cite{Cucker,Claudinei2}. Therefore, we have the next theorem:
\begin{thm}\label{symmetric}
  Let $T$ be a bounded operator given by $(\ref{T2})$ and $K$ its integral kernel given by (\ref{Kernel}). The following conditions are equivalent:
   \begin{enumerate}
       \item The kernel $K$ is symmetric.
       \item The matrix symbol $\sigma_T (\xi)$ is a hermitian matrix, for all $[\xi]\in\hat{G}$.
       \item The operator $T$ is self-adjoint.
   \end{enumerate}
\end{thm}

\pf
The proof follows by the Lemma \ref{Lm1}, \cite[Proposition 4.6]{Cucker} and definition of a self-adjoint operator in $L^2 (G)$.
\eop

\

For the positive definiteness of the kernel $K$, we state the following theorem.
\begin{thm}\label{pd}
    Let $T$ be a bounded operator of trace class given by $(\ref{T2})$ and $K$ its integral kernel given by (\ref{Kernel}). The following conditions are equivalent:
   \begin{enumerate}
    \item The kernel $K$ is positive definite.
     \item The operator $T$ is positive definite.
       \item The matrix symbol $\sigma_T (\xi)$ is positive definite for all $[\xi]\in\hat{G}$.
      
   \end{enumerate}
\end{thm}

\pf 
Let us assume that kernel $K$ is positive definite. Since $K$ is a continuous function on the compact space $G$ (Lemma \ref{Continuity}), we have that $K$ is bounded and in this way $K\in L^2 (G\times G, d\mu\times d\mu)$. Through an adaptation in the proof of Theorem 2.1 in \cite{Claudinei} for compact Lie groups endowed with a Haar measure,  we get that the bounded operator $T$ is positive definite $(1 \rightarrow 2)$. Now, if $T$ is positive definite operator and $K$ is a continuous function, then $K$ is also positive definite \cite[Theorem 2.3]{Claudinei} $(2 \rightarrow 1)$.

By the Plancherel theorem and the expression $(\ref{T})$, we have
\begin{align*}
  \langle Tf,f\rangle_{L^2 (G)} &= \langle \widehat{Tf}, \hat{f}\rangle_{\ell^2 (\hat{G})} \\
  &= \langle \widehat{f*k},\hat{f}\rangle_{\ell^2 (\hat{G})}\\
  &= \langle \sigma_T \hat{f}, \hat{f}\rangle_{\ell^2 (\hat{G})}\\
  &= \sum_{[\xi]\in\hat{G}} d_{\xi} \mathrm{Tr} [\hat{f}(\xi)^* \sigma_T (\xi) \hat{f} (\xi)].
\end{align*}
Thereby, if the matrix $\sigma_T (\xi)$ is positive definite for all $[\xi]\in\hat{G}$, then the operator $T$ is positive definite $(3\rightarrow 2)$. 

Take any $[\xi_0]\in \hat{G}$ and consider $\hat{f}\in \ell^2 (\hat{G})$ given by
\[
\hat{f}(\xi )= \left\{
\begin{array}{rll}
M(\xi_0 )\,, & \hbox{if} & \xi=\xi_0 \\
0\,, & \hbox{if} & \xi\neq \xi_0 ,
\end{array}
\right.
\]
with $M(\xi_0 )$ a matrix of size $d_{\xi_0}\times d_{\xi_0}$ such that  
\[
M(\xi_0 )_{ij} = v_i \,\delta_{j{l_0}} = \left\{
\begin{array}{rll}
v_i , & \hbox{if} & j=l_0, \\
0, & \hbox{if} & j\neq l_0,
\end{array}
\right.
\]
where $ v_i \in\mathbb{C}$ for $ 1\leq i\leq d_{\xi_0} $ and $l_0$ any fixed with $ 1\leq l_0\leq d_{\xi_0} $. Then
\[
\langle Tf_{\xi_0} ,f_{\xi_0} \rangle_{L^2 (G)}= \, d_{\xi_{0}} \mathrm{Tr} \left[M^* (\xi_0) \sigma_T (\xi_0) M(\xi_0) \right]= d_{\xi_0} \, v^* \sigma_T (\xi_0) v,
\]
where $v$ is the vector $(v_1 ,v_2 ,...,v_{d_{\xi_0}})^{\top}$. Hence, if $ T$ is positive definite then $\sigma_T (\xi)$ is a positive definite matrix $(2\rightarrow 3)$, concluding the proof.
\eop

\section{The Reproducing Kernel Hilbert Space (RKHS) associated to the kernel K}

From now on, let us consider a trace class invariant operator $T$ and its symmetric positive definite kernel $K$, given by (\ref{T2}) and (\ref{Kernel}), respectively. By Theorems \ref{symmetric} and \ref{pd}
presented in the previous section, we have that the matrix symbol $\sigma_T (\xi)$ is hermitian positive definite and, consequently, there exists a unique positive matrix $H_{\sigma_T} (\xi)$ such that $\sigma_T (\xi)= H_{\sigma_T}^2 (\xi)$ and $H_{\sigma_T} (\xi)^* =H_{\sigma_T} (\xi)$  with $[\xi]\in \hat{G}$ \cite{Koeber}. Moreover, 
\[\label{normaop}
\sup_{[\xi]\in\hat{G}} \Vert \sigma_T (\xi)\Vert_{op} 
= \sup_{[\xi]\in\hat{G}} \Vert H_{\sigma_T} (\xi) \Vert_{op}^2 ,
\]
and
\[
|\sigma_T (\xi)|=\sqrt{\sigma_T (\xi)^* \sigma_T (\xi)}=\sigma_T (\xi).
\]

Moore's theory \cite{Moore1,Moorebook}, which was followed years later by Aronszajn's work on RKHSs \cite{aron}, ensures that for a positive definite symmetric kernel $K$, there exists a unique Hilbert space \\$(\mathcal{H}_{K},\langle \,\cdot\, , \cdot \, \rangle_K)$ of functions on $G$ satisfying:
\begin{enumerate}
\item $K(\, \cdot\, ,x )\in \mathcal{H}_{K}$ for all $x\in G$;
\item (Reproducing property) $f(x)=\langle{f},{K(\cdot,x)}\rangle_{K}$, for all $x\in G$ and $f\in \mathcal{H}_{K}$.
\end{enumerate}

The following Theorem 4.2 presents the characterization of the RKHS associated to the kernel $K$.
\begin{rem}
We note that 
$$
g(x)=\sum_{[\xi]\in\hat{G}} {d_\xi } \,\mathrm{Tr} [C(\xi) \xi(x)H_{\sigma_T} (\xi)],
$$
with $C\in \ell^2 (\hat{G})$, is equivalent to  
$$
\hat{g}(\xi)= H_{\sigma_T}(\xi) C(\xi), 
$$
for all $[\xi]\in\hat{G}.$
\end{rem}
\begin{thm}\label{RKHS}
If $K:G\times G \rightarrow \mathbb{C}$ is a symmetric positive definite kernel as in (\ref{Kernel}) with $T$ a bounded continuous linear trace class operator, then 
\[
\mathcal{H}_K = \left\{g:G \rightarrow \mathbb{C} \mid
 g(x)=\sum_{[\xi]\in\hat{G}} {d_\xi } \,\mathrm{Tr} [C(\xi) \xi(x)H_{\sigma_T} (\xi)],\,\,a.e.\,\, x\in G,\, C\in \ell^2 (\hat{G})\right\}
\]
endowed with the inner product 
\[
\langle g,h\rangle_K = \sum_{[\xi]\in\hat{G}} d_{\xi} \,\mathrm{Tr} [C(\xi){B}(\xi)^*], \quad g,h \in \mathcal{H}_K,
\]
where
\[
g(x)=\sum_{[\xi]\in\hat{G}} {d_\xi } \,\mathrm{Tr} [C(\xi) \xi(x)H_{\sigma_T} (\xi)]\,\, \mbox{and}\,\, h(x)= \sum_{[\xi]\in\hat{G}} {d_\xi } \,\mathrm{Tr} [{B}(\xi) \xi(x)H_{\sigma_T}(\xi)] \,\,\mbox{for almost every} \,\,x\in G.
\]
\end{thm}

\pf Consider the set

$$
\mathcal{H} = \left\{g:G \rightarrow \mathbb{C} \mid
 g(x)=\sum_{[\xi]\in\hat{G}} {d_\xi } \,\mathrm{Tr} [C(\xi) \xi(x)H_{\sigma_T} (\xi)],\,\,a.e.\,\, x\in G,\, C\in \ell^2 (\hat{G})\right\}
$$
endowed with the inner product $\langle \cdot, \cdot \rangle_K $ as in the statement of this theorem. First we prove that if a function $g$ belongs to $\mathcal{H}$ then $g$ is an almost everywhere finite function. In fact, we prove that $g$ is a continuous function.

From the inequality (\ref{1}) ($r=1$, $p=1$ and $q=\infty$), we have that
\[
|\mathrm{Tr}[H_{\sigma_T} (\xi) C(\xi) \xi (x)] |\leq \| H_{\sigma_T} (\xi) \, C(\xi) \|_{S_1} \| \xi(x)\|_{op} = \| H_{\sigma_T} (\xi) \, C(\xi) \|_{S_1}.
\]
By the inequality (\ref{1}) again ($r=1$, $p=2$ and $q=2$), we obtain that
\[
 \| H_{\sigma_T} (\xi) \, C(\xi) \|_{S_1} \leq \| H_{\sigma_T} (\xi) \|_{S_2} \| C(\xi)\|_{S_2} = \sqrt{\mathrm{Tr}[\sigma_T (\xi)]} \, \| C(\xi)\|_{S_2} .
\]
In this way 
\[
\vert d_{\xi} \mathrm{Tr}[H_{\sigma_T} (\xi) C(\xi) \xi (x)] \vert \leq  \sqrt{d_{\xi} \, \mathrm{Tr}[\sigma_T (\xi)]} \, \sqrt{d_{\xi}}\, \| C(\xi)\|_{S_2} ,
\]
while 
\begin{align*}
\sum_{[\xi]\in\hat{G}} \sqrt{d_{\xi} \, \mathrm{Tr}[\sigma_T (\xi)]} \, \sqrt{d_{\xi}}\, \| C(\xi)\|_{S_2} &\leq \left( \sum_{[\xi]\in\hat{G}} d_{\xi} \mathrm{Tr} [\sigma_T (\xi)]\right)^{1/2} \left( \sum_{[\xi]\in\hat{G}} d_{\xi} \,\| C(\xi)\|_{S_2}^2\right)^{1/2}\\
& = \left( \sum_{[\xi]\in\hat{G}} d_{\xi} \mathrm{Tr}[\sigma_T (\xi)]\right)^{1/2} \| C\|_{\ell^2 (\hat{G})}\\
&< \infty.
\end{align*}
By the M-Test of Weierstrass, 
\[
g(x)=\sum_{[\xi]\in\hat{G}} {d_\xi } \,\mathrm{Tr} [C(\xi) \xi(x)H_{\sigma_T} (\xi)]
\]
converges uniformly in $G$ whenever $C\in \ell^2 (\hat{G})$. Therefore $g(x)$ is a continuous function on a compact set $G$ and consequently, a finite function almost everywhere.

Now, we show that the properties 1. and 2. are fulfilled. 

We write 
\[
C(\xi,y)\doteq \left(\xi(y) H_{\sigma_T}(\xi) \right)^* ,\quad y\in G, \,[\xi]\in\hat{G} .
\]
Let $L_y : \hat{G}\rightarrow \mathcal{U}({H}_\xi)$ be such that $L_y (\xi)=\xi (y)$, where $\mathcal{U}({H}_\xi)$ is the unitary group of the Hilbert space ${H}_\xi$. Since $L_y \in l^\infty (\hat{G})$ and $H_{\sigma_T}\in \ell^2(\hat{G})$, we get $L_y H_{\sigma_T}\in \ell^2 (\hat{G})$ \cite[Subsection 2.1.4]{Michael2}.
Thereby 
\[
K_{y}(x)\doteq K(x,y)= \sum_{[\xi]\in\hat{G}}{d_\xi} \,\mathrm{Tr} \left[(\xi(y) H_{\sigma_T}(\xi))^* \xi(x) H_{\sigma_T} (\xi)\right]
\]
belongs to $\mathcal{H}$ for all $x \in G$ and the property 1. is satisfied.

In order to show 2., for every $g\in \mathcal{H}$, we write
\[
g(x)=\sum_{[\xi]\in\hat{G}} {d_\xi } \,\mathrm{Tr} [C(\xi) \xi(x)H_{\sigma_T} (\xi)], \quad a.e. \,\,x\in G,
\]
and from the definition of the inner product in $\mathcal{H}$ we see that
\begin{align*}
\langle g, K_y \rangle_K &= \left\langle \sum_{[\xi]\in\hat{G}} {d_\xi} \,\mathrm{Tr} [C(\xi)  \xi(x) H_{\sigma_T}(\xi)], \sum_{[\xi]\in\hat{G}} {d_\xi} \,\mathrm{Tr} [ (\xi(y) H_{\sigma_T}(\xi))^*  \xi(x) H_{\sigma_T}(\xi)] \right\rangle_K \\ 
& =\sum_{[\xi]\in\hat{G}} {d_\xi} \,\mathrm{Tr} [C(\xi)  \xi(y) H_{\sigma_T}(\xi)]\\
&=g(y),
\end{align*}
for all $y\in G$ and $g\in \mathcal{H}$.

It is not difficult to see that $(\mathcal{H}, \langle \cdot, \cdot \rangle_K)$ is a Hilbert space and due to the uniqueness of the RKHS induced by $K$ we obtain $\mathcal{H}=\mathcal{H}_K$ and the statement of the theorem follows.
\eop

\section{Covering numbers}
We present general properties of the covering numbers we are going to use in the following subsections. Consider Banach spaces $(X, \|\cdot\|_X)$ and $(Y,\|\cdot\|_Y)$. For $\epsilon>0 $, if $B_X$ and $B_Y$ are the unit balls in $X$ and $Y$, respectively, then the \emph{covering numbers} of an operator $L: X \longrightarrow Y$ are given by
\[
 \mathcal{C} (\epsilon, L)
\doteq \mathcal{C}(\epsilon,L(B_X)) = \textrm{min} \left\{ n\in \mathbb{N} \mid \exists \,y_1, y_2,...,y_n \in Y \,\,\, \mbox{s.t.}\,\,\, L (B_X) \subset \bigcup_{j=1}^{n} (y_j + \epsilon B_Y ) \right\}.
\]
The following properties can be found in \cite{Gonzalez2, Kuhn}. If $S,\, L: X\rightarrow Y$ and $R: Z \rightarrow X$ are linear operators on Banach spaces, then for any $\epsilon, \delta>0$ we have the following properties:
\begin{enumerate}
\item[1.] $\mathcal {C} (\epsilon + \delta, L+S) \leq \mathcal{C} (\epsilon, L) \hspace{0.1cm} \mathcal{C}(\delta, S)$;
 \item[2.] $\mathcal{C}(\epsilon \delta,LR)\leq \mathcal{C}(\epsilon, L)\,\mathcal{C}(\delta, R)$;
\item[3.] If $n\doteq\mbox{rank}(L) <\infty$, then $\mathcal{C}(\epsilon, L)\leq \left( 1+2\| L\|_{op}/\epsilon \right)^ n$;
\item[4.] If $\| L\|_{op} \leq \epsilon $, then $\mathcal{C}(\epsilon, L)=1$.
\end{enumerate}
For $X$ and $Y$ two $n$-dimensional Hilbert spaces and an operator $L:X\rightarrow Y$ we have the estimate
\begin{equation}\label{CN1}
\sqrt{\det (L^* L)} \left(\frac{1}{\epsilon}\right) ^n \leq \mathcal{C}(\epsilon,L), \quad \epsilon>0.
\end{equation}
\subsection{Operator norms}\label{Opnorms} 
Properties 3. and 4. for covering numbers highlight the importance of calculating the operator norms that we present in this subsection.

If $T \in S_{1}$ is a positive trace class invariant operator on $L^{2}(G)$ as in (\ref{T2}), and $\mathfrak{F}: L^{2}(G) \rightarrow l^{2}(\hat{G})$ the group Fourier transform as in (\ref{formula2}), then
$$
\mathfrak{F} \circ T \circ \mathfrak{F}^{-1}=\bigoplus_{[\xi] \in \hat{G}} \sigma_{T}(\xi),
$$
with $\sigma_{T}$ the field of positive Hermitian operators on $\mathbb{C}^{d_{\xi}}$. 
Let the linear bounded operator $Q_{T}:l^{2}(\hat{G}) \rightarrow C(G)$ be defined by
\begin{equation}\label{QT}
Q_{T} C \doteq \mathfrak{F}^{-1}\left[H_{\sigma_{T}} C\right],
\end{equation}
for all $C \in l^{2}(\hat{G})$, where $H_{\sigma_{T}}$ is as in the Section 4,
and let $Q_{K}: l^{2}(\hat{G}) \rightarrow\left(\mathcal{H}_{K},\langle\cdot, \cdot\rangle_{K}\right)$ be an operator satisfying 
$$
Q_{K} C=Q_{T} C,
$$
for all $C \in l^{2}(\hat{G}) $. Thereby $\mathcal{H}_{K}=Q_{K}\left(l^{2}(\hat{G})\right)$. 

Note that 
\begin{align*}
\langle Q_K C,Q_K B\rangle_{K} &= \langle  \mathfrak{F}^{-1}\left[H_{\sigma_{T}} C\right] , \mathfrak{F}^{-1}\left[H_{\sigma_{T}} B\right] \rangle_K \\
&= \sum_{[\xi]\in \hat{G}} d_{\xi} \mathrm{Tr}[C(\xi) B^{*} (\xi)]\\
&= \langle C,B\rangle_{\ell^2 (\hat{G})} ,
\end{align*}
for all $C,B \in \ell^2 (\hat{G})$, thus the operator $Q_K$ is unitary.

We are interested in establishing estimates for the covering numbers $C\left(\epsilon, I_{K}\right)$ as $\epsilon \rightarrow 0$ for the embedding operator
$$
I_{K}:\left(\mathcal{H}_{K},\langle\cdot, \cdot\rangle_{K}\right) \rightarrow C(G).
$$
From now on, we will write $Q$ for $Q_T$. It is clear, that $Q=I_{K} \circ Q_{K}$, and using the unitarity of $Q_K$, we get 
\begin{align}\label{NormIK}
C\left(\epsilon, I_{K}\right)=C\left(\epsilon, Q\right),
\end{align}
for all $\epsilon>0$. 

For every subset $A \subseteq \hat{G}$, we consider the projection 
$\mathbb{P}_A : \ell^2 (\hat{G})\rightarrow \ell^2 (A)$ given as a multiplier 
$$
\mathbb{P}_{A}=\bigoplus_{[\xi] \in A} \mathbf{1}_{\mathbb{C}^{d_\xi}} .
$$
Then, we can write 
$$
g= \mathbb{P}_A \,g+ \mathbb{P}_{A^\complement} \,g, \quad g\in \ell^2 (\hat{G}) ,
$$
for any $A \subseteq\hat{G}$, and we consider the operators $Q_A \doteq Q\, \mathbb{P}_A$ and $Q_A^\complement \doteq Q\, \mathbb{P}_{A^\complement}$. Notice that $Q_{\hat{G}}=Q$.

Let $\mathcal{P}_{\sharp}(\hat{G})$ be the set of all finite subsets of $\hat{G}$. If $A\in \mathcal{P}_{\sharp}(\hat{G})$, then the operators $Q_{A}$ are finite rank operators, with 
\begin{equation}\label{Ulambda}
\mathrm{rank}(Q_A)\leq \mathrm{rank}(\mathbb{P}_A) = \sum_{[\xi] \in A} {d_\xi}^2 ,
\end{equation}
and
$$
Q-Q_{A}=Q_{A^\complement} .
$$
In particular,
$$
\left\|Q-Q_{A}\right\|_{op}=\left\|Q_A^{\complement}\right\|_{op} .
$$
More precisely, the operator norms of $Q$, $Q_A$ and $ Q_A^\complement$ are given in the following result.
\begin{lem} \label{Norms}
Let $K:G\times G \rightarrow \mathbb{C}$ a positive definite symmetric kernel as in (\ref{Kernel}) with $T$ a bounded continuous linear operator of trace class. The operators 
$$
Q : \ell^2 (\hat{G}) \longrightarrow C(G),\quad \mbox{with} \quad Q(C) \doteq \mathfrak{F}^{-1}\left[H_{\sigma_{T}} C\right],\,\, {C\in \ell^2 (\hat{G}}),
$$
$$\quad Q_A= Q\, \mathbb{P}_A, \quad\mbox{and} \quad Q_A^\complement= Q\,\mathbb{P}_{A^\complement} ,
$$
with $A\in\mathcal{P}_{\sharp}(\hat{G})$ and $C(G)$ the space of continuous functions with the supremum norm $\Vert \cdot\Vert_\infty$, satisfy respectively,
$$
\| Q \|^2_{op} = \sum_{[\xi]\in\hat{G}} d_\xi \,\mathrm{Tr}[\sigma_{T} (\xi)], \quad \|Q_A \|^2_{op} = \sum_{[\xi]\in A} d_\xi \,\mathrm{Tr}[\sigma_{T} (\xi)]
, \quad \mbox{and} \quad
\|Q_A^\complement\|^2_{op} = \sum_{[\xi]\in A^\complement} d_\xi \,\mathrm{Tr}[\sigma_{T} (\xi)].
$$
\end{lem}

\pf From Lemma 10.3.28 in \cite{Michael1} and inequality (\ref{1}) ($r=2$, $p=2$, and $q=\infty$) we have
$$
\left| \sum_{[\xi]\in\hat{G}} d_{\xi} \mathrm{Tr} [\xi(x)H_{\sigma_T} (\xi) C(\xi)] \right|\leq \sum_{[\xi]\in\hat{G}} d_{\xi} \Vert H_\sigma (\xi)\Vert_{S^2} \Vert C(\xi)\Vert_{S^2},
$$
for all $x\in G$, and from Cauchy--Schwarz inequality we obtain
$$
\sum_{[\xi]\in\hat{G}} d_{\xi} \Vert H_\sigma (\xi)\Vert_{S^2} \Vert C(\xi)\Vert_{S^2} \leq \left( \sum_{[\xi]\in\hat{G}} d_\xi \Vert H_{\sigma} (\xi)\Vert^2_{S^2} \right)^{1/2} \left( \sum_{[\xi]\in\hat{G}} d_\xi \Vert C(\xi)\Vert^2_{S^2} \right)^{1/2} =\Vert H_{\sigma}\Vert _{\ell^2(\hat{G})} \Vert C\Vert _{\ell^2(\hat{G})} ,
$$
for all $C\in \ell^2 (\hat{G})$, then
$$
\Vert Q(C)\Vert_{\sup} =\Vert (\mathfrak{F}^{-1}\left[H_{\sigma_{T}} C\right]) \Vert_{\sup} \leq \Vert H_{\sigma_T}\Vert _{\ell^2(\hat{G})} \Vert C\Vert _{\ell^2(\hat{G})},
$$
for all $C\in \ell^2 (\hat{G})$. Thereby,
$$
\left\|Q\right\|_{op} \leq \left\|H_{\sigma_{T}}\right\|_{\ell^2 (\hat{G})} .
$$
If $C=H_{\sigma_{T}}$ and $x=e$ (the unit element of $G$), then 
$$
(Q(C))(x)= (\mathfrak{F}^{-1}\left[{H_{\sigma_{T}}}^2\right])(e)= \sum_{[\xi]\in \hat{G}} d_{\xi} \mathrm{Tr}[\sigma_T (\xi)]= \Vert H_{\sigma_T}\Vert_{\ell^2 (\hat{G})}^2 ,
$$
and the upper limit is attained. Therefore,
$$
\left\|Q\right\|_{op}=\|H_{\sigma_T}\|_{\ell^2(\hat{G})}=\left(\sum_{[\xi]\in\hat{G}} d_{\xi} \mathrm{Tr} [\sigma_T (\xi)]\right)^{1/2} .
$$
The remaining norm equalities follow by analogous calculations. 
\eop

Since by \cite[Theorem 5.1]{Delgado2} we have $\Vert H_{\sigma_T}\Vert _{\ell^2(\hat{G})} ^2  = \sum_{[\xi]\in\hat{G}} d_{\xi} \mathrm{Tr} [\sigma_T (\xi)] =\mathrm{Tr} \,T$ , we get 
$$
\left\|H_{\sigma_{T}}\right\|_{\ell^2 (\hat{G})} 
=\sqrt{\|T\|_{S^1}} ,
$$ 
so that
$$
\left\|Q\right\|_{op} =\sqrt{\|T\|_{S^1}} .
$$
Rewriting $\epsilon= \epsilon \,\Vert Q \Vert_{op} / \Vert Q \Vert_{op}$ and $Q = Q  \Vert Q \Vert_{op}/ \Vert Q \Vert_{op}$, by the properties $2.$ and $4.$ for covering numbers, we have 
\begin{align}\label{DivNormQ}
C\left(\epsilon, Q\right)=C\left(\frac{\epsilon}{\,\,\,\,\left\|Q\right\|_{op}}, \frac{Q}{\,\,\,\,\left\|Q\right\|_{op}}\right),
\end{align}
and we can assume for simplicity that $\left\|Q\right\|_{op}=\|T\|_{S^1}=1$.

\subsection{Operator rank estimates}
In the previous subsection, finite rank operators $Q_A$ for some $A\in \mathcal{P}_{\sharp}(\hat{G})$ were considered. In this subsection, we consider a smaller family of finite sets $A$ and we estimate the rank for the associated operator.

Let $\mathcal{L}_G $ be the Laplacian operator on $G$ (or the Casimir element of the universal enveloping algebra).
For each $[\xi] \in \hat{G}$, the matrix elements of $\xi$ are the eigenfunctions for $\mathcal{L}_G $ with eigenvalues denoted by $ -\lambda_{[\xi]}$ \cite[Theorem 10.3.13]{Michael1}. In other words,
\[
-\mathcal{L}_G \xi_{ij} (x)=\lambda_{[\xi]} \xi_{ij} (x),\quad x\in G,
\]
where $\lambda_{[\xi]}\geq 0$ and $1\leq i,j\leq d_{\xi}$. Moreover,
 $\langle\xi\rangle \doteq (1+\lambda_{[\xi]})^{1/2}$ is an eigenvalue of the first-order elliptic pseudo-differential operator $\left(I-\mathcal{L}_G\right)^{1 / 2}$.

Instead of optimising over all possible finite subsets $A \in \mathcal{P}_{\sharp}(\hat{G})$, we optimise over the family of finite sets $\left\{A_{\lambda}\right\}_{\lambda \in(1,+\infty)}$ given by
\begin{equation}\label{Alambda}
A_{\lambda} \doteq\{[\xi] \in \hat{G} \mid \quad\langle\xi\rangle \leq \lambda\},
\end{equation}
for all $\lambda \in(1,+\infty)$.
It is clear that the family $\left\{A_{\lambda}\right\}_{\lambda \in(1,+\infty)}$ is increasing with $\lambda$,
$$
A_{\lambda} \subset A_{\mu},\quad \lambda \leq \mu ,
$$
for all $\lambda, \mu \in(0,+\infty)$. In this way, Lemma \ref{Norms} holds for finite sets $\left\{A_{\lambda}\right\}_{\lambda \in(1,+\infty)}$ as in (\ref{Alambda}) and for operators
$$
Q_{A_\lambda}= Q\, \mathbb{P}_{A_\lambda}, \quad\mbox{and} \quad Q_{A_\lambda}^\complement= Q\, \mathbb{P}_{A_\lambda^\complement} ,
$$
we have
$$
\|Q_{A_\lambda} \|^2_{op} = \sum_{[\xi]\in {A_\lambda}} d_\xi \,\mathrm{Tr}[\sigma_{T} (\xi)]= \sum_{\langle\xi\rangle \leq \lambda} d_\xi \,\mathrm{Tr}[\sigma_{T} (\xi)] ,
$$
and 
$$
\|Q_{A_\lambda}^\complement\|^2_{op} = \sum_{[\xi]\in A_\lambda^\complement} d_\xi \,\mathrm{Tr}[\sigma_{T} (\xi)]= \sum_{\langle\xi\rangle > \lambda} d_\xi \,\mathrm{Tr}[\sigma_{T} (\xi)].
$$
Moreover, we can write the formula (\ref{Ulambda}) as
$$
\mathrm{rank}(Q_{A_\lambda})\leq\sum_{\langle\xi\rangle \leq \lambda} {d_\xi}^2 .
$$
\begin{rem} We record two properties that will be useful on some occasions:
\begin{equation*}
\sum_{\langle\xi\rangle \leq \lambda} {d_{\xi}}^2 \langle \xi\rangle^{\alpha n} \asymp \lambda^{(\alpha+1)n} \quad \mbox{for}\,\, \alpha> -1,
\end{equation*}

\begin{equation*}
\sum_{\langle\xi\rangle \geq \lambda} {d_{\xi}}^2 \langle \xi\rangle^{\alpha n} \asymp \lambda^{(\alpha+1)n} \quad \mbox{for}\,\, \alpha< -1, 
\end{equation*}
as $\lambda \rightarrow \infty$. The formulas mentioned above, found in \cite{Delgado4}, are originally derived in the asymptotic sense as $\lambda\rightarrow+\infty$. However, it is straightforward to observe that these formulas also hold uniformly for $\lambda\in(1,+\infty)$. So that, there exist real positive constants $C_{\alpha,n}$, $c_{\alpha,n}$, $W_{\alpha,n} $ and $w_{\alpha,n}$ such that
\begin{equation}\label{assym1}
c_{\alpha,n} \lambda^{(\alpha+1)n}\leq\sum_{\langle\xi\rangle \leq \lambda} {d_{\xi}}^2 \langle \xi\rangle^{\alpha n} \leq C_{\alpha,n} \lambda^{(\alpha+1)n}, \quad \mbox{for}\,\, \alpha> -1,
\end{equation}

\begin{equation}\label{assym2}
w_{\alpha,n} \lambda^{(\alpha+1)n} \leq \sum_{\langle\xi\rangle \geq \lambda} {d_{\xi}}^2 \langle \xi\rangle^{\alpha n} \leq W_{\alpha,n} \lambda^{(\alpha+1)n}, \quad \mbox{for}\,\, \alpha< -1, 
\end{equation}
for all $\lambda\in(1,+\infty)$.
\end{rem}
Note that by the formula (\ref{assym1}) with $\alpha=0$, one obtains that 
\begin{equation}\label{dim}
\mathrm{rank}(Q_{A_\lambda}) \leq C_{0,n}\lambda^n ,\quad \lambda\in(1,+\infty).
\end{equation}

\subsection{Estimates related to kernels with symbol with trace and determinant of finite order} 
In this subsection, we consider $T$ a positive trace class invariant operator on $L^2 (G)$ as in (\ref{T2}). In order to find estimates for covering numbers of the operator $Q$, and as a consequence of the $I_K$ operator, we define an order notion to the trace and the determinant
of the symbol $ \sigma_T$ as follows. 
\begin{defn}
Let $\beta > 0$. The trace of the symbol $\sigma_T$ has order less than or equal to $\beta$ if there exist real positive constants $b_T$ such that 
\[
\mathrm{Tr}[\sigma_T (\xi)] \leq b_T \,d_{\xi} \, \langle \xi \rangle^{-\beta}
\]
for all $[\xi]\in\hat{G}$, respectively.
\end{defn}
We note that if 
the trace of the symbol $\sigma_T$ has order less than or equal to $\beta$ then
\begin{equation}\label{beta}
  \sum_{[\xi]\in\hat{G}} d_{\xi} \mathrm{Tr}[\sigma_T (\xi)] \leq b_T \sum_{[\xi]\in\hat{G}} {d_{\xi}}^2 \langle \xi\rangle^{-\beta},
\end{equation}
and if we assume $\beta> n$, then by \cite[Lemma 3.8]{Delgado2} we obtain  
$$
\sum_{[\xi]\in\hat{G}} {d_{\xi}}^2 \langle \xi\rangle^{-\beta} < \infty,
$$
and $T$ is a trace class operator. 
\begin{defn}
 Let $ \gamma > 0$. The determinant of the symbol $\sigma_T$ has order greater than or equal to $\gamma$ 
if there exist real positive constants $c_T$ and $\omega_T $ such that 
$$
\left(\det[\sigma_T (\xi)]\right)^{\frac{1}{d_{\xi}}} \geq c_T \,e^{-2 \omega_T\langle \xi \rangle^{\gamma}},
$$
for all $[\xi]\in\hat{G}$, respectively.   
\end{defn}

\begin{thm}\label{thmuppergen}  Let
 $K$ be a positive definite symmetric kernel on a n-dimensional compact Lie group $G$, represented by the series expansion (\ref{Kernel}).\ If there exists $\beta >n$ such that the trace of the symbol $\sigma_T$ has order less than or equal to $\beta$ 
 then there exist positive real constants $C_n$, $b_T$ and $\kappa_\beta$ such that 
$$
\ln C\left(\epsilon, I_{K}\right)
\leq \frac{C_{n}\left(4 b_{T} \kappa_{\beta}\|T\|_{S^1}\right)^{\frac{n}{\beta-n}}}{\epsilon^{\frac{2 n}{\beta-n}}} \ln \left[1+\frac{4 \sqrt{\|T\|_{S^1}}}{\epsilon}\right],
$$
for all $\epsilon \in\left(0, \frac{\sqrt{\|T\|_{S^1}}}{\sqrt{3}}\right)$.
\end{thm}

\pf We consider approximations $Q_{A_\lambda}\doteq Q\, \mathbb{P}_{A_\lambda}$ with $A_\lambda$ as in formula (\ref{Alambda}), to the operator $Q$, such that $\Vert Q -Q_{A_\lambda}\Vert \leq \delta$, for all $ \lambda \in (1, +\infty)$. In other words, 
$$
\sum_{[\xi] \in A_{\lambda}^\complement} d_{\xi} \operatorname{Tr}\left[\sigma_{T}(\xi)\right] \leq \delta^2 , \quad  \lambda \in (1, +\infty),\quad \delta >0 .
$$
By the property 1. of covering numbers we write
$$
\mathcal{C} (\epsilon, Q)\leq \mathcal{C}(\epsilon-\delta, Q_{A_\lambda}) \mathcal{C} (\delta, Q -  Q_{A_\lambda}), \quad 0<\delta<\epsilon,
$$
and by the properties 4. and 3. of covering numbers we obtain
$$
\mathcal{C} (\epsilon, Q)\leq \mathcal{C}(\epsilon-\delta, Q_{A_\lambda})\leq \left(1+\frac{2}{\epsilon-\delta} \sqrt{\sum_{[\xi] \in A_\lambda} d_{\xi} \operatorname{Tr}\left[\sigma_{T}(\xi)\right]}\right)^{\mathrm{rank} \left(Q_{A_\lambda}\right)} , \quad 0<\delta<\epsilon.
$$
Denoting
$$
\delta_{\lambda} \doteq \sqrt{\sum_{[\xi] \in A_{\lambda}^{\complement}} d_{\xi} \operatorname{Tr}\left[\sigma_{T}(\xi)\right]}, \quad \lambda \in(1,+\infty),
$$
we see that $\delta_\lambda\leq\delta <\epsilon $. Thereby, if we consider $\Vert Q\Vert =1$ as at the end of Subsection \ref{Opnorms}, we have
$$
\mathcal{C} (\epsilon, Q)\leq\left(1+\frac{2 \sqrt{1-\delta_{\lambda}^{2}}}{\epsilon-\delta_{\lambda}}\right)^{\mathrm{rank} \left(Q_{A_\lambda}\right)} \quad, \quad \lambda \in(1,+\infty), \quad \delta_{\lambda}<\epsilon <1.
$$
Assuming $\epsilon<\frac{1}{\sqrt{3}}$, the function $F:(0, \epsilon) \rightarrow(0,+\infty)$ given by
$$
F(\delta) \doteq \frac{\sqrt{1-\delta^{2}}}{\epsilon-\delta},
$$
for all $ \delta \in(0, \epsilon)$
is strictly increasing. If we assume that $\beta >n$ such that the trace of the symbol $\sigma_T$ has order less than or equal to $\beta$, then
$$
\operatorname{Tr}\left[\sigma_{T}(\xi)\right] \leq b_{T} \,d_{\xi}\langle\xi\rangle^{-\beta}, \quad \beta \in(n,+\infty),\quad [\xi]\in \hat{G},
$$
for some constant $b_T \in(0,+\infty)$. Then by ($\ref{assym2}$) $(\alpha=-\beta/n)$ we obtain
$$
\delta_{\lambda}^{2}=\sum_{\langle\xi\rangle>\lambda} d_{\xi} \operatorname{Tr}\left[\sigma_{T}(\xi)\right] \leq b_{T} \sum_{\langle\xi\rangle > \lambda} d_{\xi}^{2}\langle\xi\rangle^{-\beta} \leq b_{T} \kappa_{\beta} \lambda^{n-\beta},
$$
for all $\lambda \in(1,+\infty)$ and for some constant $\kappa_{\beta} \in(0,+\infty)$. It follows that
$$
F\left(\delta_{\lambda}\right) \leq \frac{\sqrt{1-b_{T} \kappa_{\beta} \lambda^{n-\beta}}}{\epsilon-\sqrt{b_{T} \kappa_{\beta}} \lambda^{\frac{n-\beta}{2}}}, 
$$
for all $\lambda \in(1,+\infty)$. Note that the restriction $\delta_{\lambda}<\epsilon$ is satisfied if we assume that
$$
\lambda>\left(\frac{b_{T} \kappa_{\beta}}{\epsilon^{2}}\right)^{\frac{1}{\beta-n}}.
$$
On the other hand, by (\ref{dim})
$$
\mathrm{rank}(Q_{A_\lambda})\leq C_{n} \lambda^{n},
$$
for all $\lambda \in(1,+\infty)$ and for some constant $C_{n} \in(0,+\infty)$. Putting all together we arrive at
$$
\mathcal{C} (\epsilon, Q)\leq\left(1+\frac{2 \sqrt{1-\delta_{\lambda}^{2}}}{\epsilon-\delta_{\lambda}}\right)^{\mathrm{rank} \left(Q_{A_\lambda}\right)}\leq \left(1+\frac{2 \sqrt{1-b_{T} \kappa_{\beta} \lambda^{n-\beta}}}{\epsilon-\sqrt{b_{T} \kappa_{\beta}} \lambda^{\frac{n-\beta}{2}}}\right)^{C_{n} \lambda^{n}},
$$
for all $\lambda \in(1,+\infty)$ and we need to minimise the function $H_\epsilon :(1,+\infty)\rightarrow \mathbb{R}$ given by
$$
H_{\epsilon}(\lambda) \doteq\left(1+\frac{2 \sqrt{1-b_{T} \kappa_{\beta} \lambda^{n-\beta}}}{\epsilon-\sqrt{b_{T} \kappa_{\beta}} \lambda^{\frac{n-\beta}{2}}}\right)^{C_{n} \lambda^{n}},
$$
defined for all $\lambda \in(1,+\infty)$,
on the interval
$$
\lambda \in\left(\left(\frac{b_{T} \kappa_{\beta}}{\epsilon^{2}}\right)^{\frac{1}{\beta-n}},+\infty\right) .
$$
If we call
$$
\lambda_{\epsilon} \doteq\left(\frac{4 b_{T} \kappa_{\beta}}{\epsilon^{2}}\right)^{\frac{1}{\beta-n}} ,
$$
then
$$
\min _{\lambda \in(1,+\infty)} \ln H_{\epsilon}(\lambda)\leq \ln H_\epsilon (\lambda_\epsilon)<C_{n} \lambda_{\epsilon}^{n} \ln \left[1+\frac{2}{\epsilon-\sqrt{b_{T} \kappa_{\beta}} \lambda_{\epsilon}^{\frac{n-\beta}{2}}}\right]=\frac{C_{n}\left(4 b_{T} \kappa_{\beta}\right)^{\frac{n}{\beta-n}}}{\epsilon^{\frac{2 n}{\beta-n}}} \ln \left[1+\frac{4}{\epsilon}\right] .
$$
From (\ref{NormIK}) and (\ref{DivNormQ}), we get that $\ln C\left(\epsilon, I_{K}\right)=\ln C\left(\frac{\epsilon}{\sqrt{\|T\|_{S^1}}}, \frac{Q}{\left\|Q\right\|}\right) $, where $\| T\|_{S^1}=\sum_{[\xi]\in\hat{G}} d_{\xi} \mathrm{Tr} [\sigma_T (\xi)] $. Thus, if 
$$
\operatorname{Tr}\left[\sigma_{T}(\xi)\right] \leq b_{T}\|T\|_{S^1} d_{\xi}\langle\xi\rangle^{-\beta} 
$$
with $\beta \in(n,+\infty)$, we arrive at the upper estimate 
\begin{align}
\ln C\left(\epsilon, I_{K}\right)&=\ln C\left(\frac{\epsilon}{\sqrt{\|T\|_{S^1}}}, \frac{Q}{\left\|Q\right\|}\right) \\
&\leq \frac{C_{n}\left(4 b_{T} \kappa_{\beta}\|T\|_{S^1}\right)^{\frac{n}{\beta-n}}}{\epsilon^{\frac{2 n}{\beta-n}}} \ln \left[1+\frac{4 \sqrt{\|T\|_{S^1}}}{\epsilon}\right],
\end{align}
for all $\epsilon \in\left(0, \frac{\sqrt{\|T\|_{S^1}}}{\sqrt{3}}\right)$.
\eop

To find lower estimates for the covering numbers of the operator $I_K : \mathcal{H}_K \rightarrow C(G)$ we consider the determinant of the symbol $\sigma_T$ of order greater than or equal to $\gamma$ as follows.

\begin{thm}\label{thmlowergeogen1}  Let
$T$ be a bounded trace class invariant operator given by (\ref{T2}) and $K$ be a positive definite symmetric kernel on a n-dimensional compact Lie group $G$ represented by the series expansion (\ref{Kernel}).\ If the determinant of the symbol $\sigma_T$ has order greater than or equal to $\gamma$, then there exist positive real constants $c_{0,n}$, $\omega_T$, $a_T$ and $\mu_\gamma$ such that 
\[
\ln C\left(\epsilon, I_{K}\right)\geq\left(\frac{c_{0,n}}{\omega_T\mu_{\gamma}\left(1+\frac{\gamma}{n}\right)}\right)^{\frac{n}{\gamma}} \frac{c_{0,n}}{1+\frac{n}{\gamma}}\left(\ln \left[\frac{a_{T} \sqrt{\|T\|_{S^1}}}{\epsilon}\right]\right)^{1+\frac{n}{\gamma}},
\]
for all $\epsilon \in\left(0, a_{T} \sqrt{\|T\|_{S^1}} e^{-\frac{\omega_T \mu \gamma \left(1+\frac{\gamma}{n}\right)}{c_{0,n}}}\right)$.
\end{thm}

\pf Let $A_\lambda$ be as in (\ref{Alambda}), we consider the composition operator as follows
\begin{equation*}\label{L_A}
L_{A_\lambda}:l^{2}(A_\lambda) \stackrel{\jmath_{\lambda}}{\hookrightarrow} l^{2}(\hat{G}) \xrightarrow{Q} C(G) \stackrel{\jmath}{\hookrightarrow}L^{2}(G) \xrightarrow{P_\lambda} F_\lambda,
\end{equation*}
where $\jmath_{\lambda}$ and $\jmath$ stand for the embeddings given by $\ell^2 (A_{\lambda})\hookrightarrow \ell^2 (\hat{G})$ and $C(G)\hookrightarrow L^2 (G)$, respectively, and $P_\lambda$ is the orthogonal projection of $L^2 (G)$ on $F_\lambda\doteq \jmath \,Q\,\jmath_{\lambda} (\ell^2 (A_\lambda))$. It is clear that $\Vert P_\lambda\Vert_{op}=\Vert \jmath\Vert_{op}=\Vert \jmath_{A_\lambda}\Vert_{op}=1$ and the operator 
$$
L_{A_\lambda}= P_\lambda \circ\jmath \circ Q \circ \jmath_{\lambda}:\ell^2 (A_{\lambda}) \rightarrow F_\lambda
$$
is bijective. By properties 2. and 4. for covering numbers, 
$$
C\left(\epsilon, L_{A_\lambda}\right) \leq  C(1, P_\lambda) C(1, \jmath) C\left(\epsilon, Q\right) C\left(1, \jmath_{{\lambda}}\right)=C\left(\epsilon, Q\right),
$$
for all $\epsilon\in(0, +\infty)$. Since $\dim \ell^2 (A_\lambda)= \dim (F_\lambda)$, we can use the lower bound for the covering number in inequality (\ref{CN1}), so that
$$
\frac{\sqrt{\operatorname{det}\left(L_{A_\lambda}^{*} L_{A_\lambda}\right)}}{\epsilon^{\dim \ell^2 (A_\lambda) }} \leq C\left(\epsilon, L_{A_\lambda}\right), 
$$
for all $\lambda\in(1, +\infty)$ and $\epsilon\in(0,+\infty) $. \\By a direct calculation, 
$$
L_{A_\lambda}^* L_{A_\lambda} C=\sigma_T C,
$$
for all $C\in \ell^2 (A_\lambda)$ with $\lambda\in(1,+\infty)$, thereby
$$
\operatorname{det}\left(L_{A_{\lambda}}^{*} L_{A_\lambda}\right)=\prod_{[\xi] \in A_\lambda}\left(\operatorname{det} \sigma_{T}(\xi)\right)^{d_{\xi}} ,
$$
for all $\lambda\in(1,+\infty)$, and 
\[
\frac{\prod\limits_{\langle\xi\rangle \leq\lambda}\left(\operatorname{det} \sigma_{T}(\xi)\right)^{\frac{d_{\xi}}{2}}}{\epsilon^{\dim \ell^2 (A_\lambda) }} \leq C\left(\epsilon, L_{A_\lambda}\right),
\]
for all $\lambda\in(1, +\infty)$ and $\epsilon\in(0,+\infty) $.\\
By the formula (\ref{assym1}) with $\alpha=0$, 
$$
\dim \ell^2 (A_\lambda)= \sum_{\langle\xi\rangle \leq \lambda} d_{\xi}^{2} \geq c_{0,n} \lambda^{n},
$$
for all $\lambda\in(1, +\infty)$ and some constant $c_{0,n}\in (0, +\infty)$. If we assume that the determinant of the symbol $\sigma_T$ has order greater than or equal to $\gamma$, then
$$
\left(\det\sigma_T (\xi)\right)^{\frac{1}{d_\xi}}\geq a_T^2 e^{-2 \omega_T \langle \xi \rangle^{\gamma}}, \quad \gamma\in(0,+\infty),
$$
for some constants $a_T , \omega_T \in(0,+\infty)$, and by (\ref{assym1}) $(\text{if}\,\, \alpha= \gamma/n)$ follows that
\begin{align*}
\ln \left[\prod_{\langle\xi\rangle \leq \lambda}\left(\operatorname{det} \sigma_{T}(\xi)\right)^{\frac{d_{\xi}}{2}}\right] &=\sum_{\langle\xi\rangle \leq \lambda} d_{\xi}^{2} \ln \left[\left(\operatorname{det} \sigma_{T}(\xi)\right)^{\frac{1}{2 d_{\xi}}}\right]\\
&\geq \ln a_{T} \sum_{\langle\xi\rangle 
\leq \lambda} d_{\xi}^{2}- \omega_T \sum_{\langle\xi\rangle \leq \lambda} d_{\xi}^{2}\langle\xi\rangle^{\gamma} \\
& \geq c_{0,n} (\ln a_{T}) \lambda^{n}-\omega_T\mu_{\gamma} \lambda^{n+\gamma}, 
\end{align*}
for all $\lambda \in(1,+\infty)$ and some constant $\mu_\gamma\in(0,+\infty)$.
Thereby, 
$$
\ln\left[\frac{\prod_{\langle\xi\rangle \leq\lambda}\left(\operatorname{det} \sigma_{T}(\xi)\right)^{\frac{d_{\xi}}{2}}}{\epsilon^{\dim \ell^2 (A_\lambda) }} \right]\geq -\omega_T \mu_{\gamma} \lambda^{n+\gamma}+c_{0,n} \lambda^{n} \ln \left[\frac{a_{T}}{\epsilon}\right],
$$
for all $\lambda\in(1,+\infty)$, and we need to maximise the function $G_\epsilon:(1, +\infty)\rightarrow \mathbb{R}$ given by
$$
G_{\epsilon}(\lambda) \doteq-\omega_T \mu_{\gamma} \lambda^{n+\gamma}+c_{0,n} \lambda^{n} \ln \left[\frac{a_{T}}{\epsilon}\right],\quad  \lambda \in(1,+\infty).
$$
Then, for
$$
\epsilon<a_{T} e^{-\frac{\mu_{\gamma} \left(1+\frac{\gamma}{n}\right)}{c_{0,n}}},
$$
we find by maximisation that
$$
\max _{\lambda \in(1,+\infty)} G_{\epsilon}(\lambda)=\left(\frac{c_{0,n}}{\omega_T\mu_{\gamma}\left(1+\frac{\gamma}{n}\right)}\right)^{\frac{n}{\gamma}} \frac{c_{0,n}}{1+\frac{n}{\gamma}}\left(\ln \left[\frac{a_{T}}{\epsilon}\right]\right)^{1+\frac{n}{\gamma}} .
$$
Thus, if 
$$
\left(\operatorname{det} \sigma_{T}(\xi)\right)^{\frac{1}{d_{\xi}}} \geq a_{T}^{2}\|T\|_{S^1} e^{-2 \omega_T\langle\xi\rangle^{\gamma}},
$$
with $\omega_T, a_{T}, \gamma \in(0,+\infty)$, we finally establish the lower estimate
\begin{align*}
\ln C\left(\epsilon, I_{K}\right)&=\ln C\left(\epsilon, Q\right)=\ln C\left(\frac{\epsilon}{\sqrt{\|T\|_{1}}}, \frac{Q}{\,\,\,\,\left\|Q\right\|_{op}}\right) \\
&\geq\left(\frac{c_{0,n}}{\omega_T\mu_{\gamma}\left(1+\frac{\gamma}{n}\right)}\right)^{\frac{n}{\gamma}} \frac{c_{0,n}}{1+\frac{n}{\gamma}}\left(\ln \left[\frac{a_{T} \sqrt{\|T\|_{1}}}{\epsilon}\right]\right)^{1+\frac{n}{\gamma}},
\end{align*}
for all $\epsilon \in\left(0, a_{T} \sqrt{\|T\|_{S^1}} e^{-\frac{\omega_T \mu \gamma \left(1+\frac{\gamma}{n}\right)}{c_{0,n}}}\right)$.
\eop

\vspace{0.7cm}

\textbf{Acknowledgements.} The first and third named authors were supported by the FWO Senior Research Grant G022821N, and by
 the Methusalem programme of the Ghent University Special Research Fund (BOF) (Grant number 01M01021). The work of the second author was supported by the Higher Education Science Committee of the Republic of Armenia (Research project no. 23RL-1A027) and by the Morá Miriam Rozen Gerber fellowship from the Weizmann Institute of Science. The second author also wishes to express her sincere gratitude to Dr. Duván Cardona for valuable discussions during her visit to the Ghent Analysis and EDP Center in Belgium.

\bibliographystyle{plain}
\bibliography{references}

@article {Ruskai,
    AUTHOR = {Ruskai, M. B.},
     TITLE = {Inequalities for traces on von {N}eumann algebras},
   JOURNAL = {Comm. Math. Phys.},
  FJOURNAL = {Communications in Mathematical Physics},
    VOLUME = {26},
      YEAR = {1972},
     PAGES = {280--289},
      ISSN = {0010-3616,1432-0916},
   MRCLASS = {46L10 (81.46)},
  MRNUMBER = {312284},
MRREVIEWER = {H.\ Wakita},
       URL = {http://projecteuclid.org/euclid.cmp/1103858124},
}

@article {zhou,
    AUTHOR = {Zhou, D.X.},
     TITLE = {Capacity of reproducing kernel spaces in learning theory},
   JOURNAL = {IEEE Trans. Inform. Theory},
  FJOURNAL = {Institute of Electrical and Electronics Engineers.
              Transactions on Information Theory},
    VOLUME = {49},
      YEAR = {2003},
    NUMBER = {7},
     PAGES = {1743--1752},
      ISSN = {0018-9448},
   MRCLASS = {62G05 (46E22 68Q32 68T05)},
  MRNUMBER = {1985575},
MRREVIEWER = {Ivan K\v{r}iv\'{y}},
       DOI = {10.1109/TIT.2003.813564},
       URL = {https://doi.org/10.1109/TIT.2003.813564},
}

@article {Kuhn,
    AUTHOR = {K\"{u}hn, T.},
     TITLE = {Covering numbers of {G}aussian reproducing kernel {H}ilbert
              spaces},
   JOURNAL = {J. Complexity},
  FJOURNAL = {Journal of Complexity},
    VOLUME = {27},
      YEAR = {2011},
    NUMBER = {5},
     PAGES = {489--499},
      ISSN = {0885-064X},
   MRCLASS = {68Q32 (41A46 46E22 60G15 68T05 94A17)},
  MRNUMBER = {2805533},
MRREVIEWER = {Ingo Steinwart},
       DOI = {10.1016/j.jco.2011.01.005},
       URL = {https://doi.org/10.1016/j.jco.2011.01.005},
}

@article{Mercer1909,
  title={Functions of Positive and Negative Type, and their Connection with the Theory of Integral Equations},
  author={J. Mercer},
  journal={Philosophical Transactions of the Royal Society A},
  year={1909},
  volume={209},
  pages={415-446},
  url={https://api.semanticscholar.org/CorpusID:121070291}
}

@book{cristianini_shawe-taylor_2000, place={Cambridge}, title={An Introduction to Support Vector Machines and Other Kernel-based Learning Methods}, DOI={10.1017/CBO9780511801389}, publisher={Cambridge University Press}, author={Cristianini, N. and Shawe-Taylor, J.}, year={2000}}

@book{Vapnik1998,
  author = {Vapnik, V.N.},
  publisher = {Wiley-Interscience},
  title = {Statistical Learning Theory},
  year = 1998
}

@article {jordao-men,
    AUTHOR = {Jordão, T. and Menegatto, V. A.},
     TITLE = {Kolmogorov widths on the sphere via eigenvalue estimates for
              {H}\"{o}lderian integral operators},
   JOURNAL = {Results Math.},
  FJOURNAL = {Results in Mathematics},
    VOLUME = {74},
      YEAR = {2019},
    NUMBER = {2},
     PAGES = {Paper No. 74, 18},
      ISSN = {1422-6383,1420-9012},
   MRCLASS = {41A46 (42A16 45C05 47B34 47G10)},
  MRNUMBER = {3923505},
MRREVIEWER = {Karol\ Dziedziul},
       DOI = {10.1007/s00025-019-1000-4},
       URL = {https://doi.org/10.1007/s00025-019-1000-4},
}

@article {Gonzalez,
    AUTHOR = {Azevedo, D. and Gonzalez, K. and Jordão, T.},
     TITLE = {Sharp estimates for the covering numbers of the {W}eierstrass
              fractal kernel},
   JOURNAL = {J. Complexity},
  FJOURNAL = {Journal of Complexity},
    VOLUME = {74},
      YEAR = {2023},
     PAGES = {Paper No. 101692},
      ISSN = {0885-064X},
   MRCLASS = {47B06 (26A15 42A16 42A32 46E22)},
  MRNUMBER = {4513086},
       DOI = {10.1016/j.jco.2022.101692},
       URL = {https://doi.org/10.1016/j.jco.2022.101692},
}

@article {kolmogorov,
    AUTHOR = {Kolmogorov, A. N. and Tikhomirov, V. M.},
     TITLE = {{$\varepsilon $}-entropy and {$\varepsilon $}-capacity of sets in function spaces},
   JOURNAL = {Uspehi Mat. Nauk},
  FJOURNAL = {Akademiya Nauk SSSR i Moskovskoe Matematicheskoe Obshchestvo.
              Uspekhi Matematicheskikh Nauk},
    VOLUME = {14},
      YEAR = {1959},
    NUMBER = {2 (86)},
     PAGES = {3--86},
      ISSN = {0042-1316},
   MRCLASS = {46.00 (26.00)},
  MRNUMBER = {0112032},
MRREVIEWER = {G. G. Lorentz},
}

@article {aron,
    AUTHOR = {Aronszajn, N.},
     TITLE = {Theory of reproducing kernels},
   JOURNAL = {Trans. Amer. Math. Soc.},
  FJOURNAL = {Transactions of the American Mathematical Society},
    VOLUME = {68},
      YEAR = {1950},
     PAGES = {337--404},
      ISSN = {0002-9947},
   MRCLASS = {46.0X},
  MRNUMBER = {51437},
MRREVIEWER = {T. H. Hildebrandt},
       DOI = {10.2307/1990404},
       URL = {https://doi.org/10.2307/1990404},
}

@incollection {Quang,
    AUTHOR = {Minh, H.Q. and Niyogi, P. and Yao, Y.},
     TITLE = {Mercer's theorem, feature maps, and smoothing},
 BOOKTITLE = {Learning theory},
    SERIES = {Lecture Notes in Comput. Sci.},
    VOLUME = {4005},
     PAGES = {154--168},
 PUBLISHER = {Springer, Berlin},
      YEAR = {2006},
   MRCLASS = {68T05 (62G08)},
  MRNUMBER = {2280604},
       DOI = {10.1007/11776420\_14},
       URL = {https://doi.org/10.1007/11776420_14},
}

@book {ingobook,
    AUTHOR = {Steinwart, I. and Christmann, A.},
     TITLE = {Support vector machines},
    SERIES = {Information Science and Statistics},
 PUBLISHER = {Springer, New York},
      YEAR = {2008},
     PAGES = {xvi+601},
      ISBN = {978-0-387-77241-7},
   MRCLASS = {62-02 (60B11 60G15 62Gxx 62H30 62J02 68T05 68T20)},
  MRNUMBER = {2450103},
MRREVIEWER = {Gilles Blanchard},
}

@article {smale,
    AUTHOR = {Cucker, F. and Smale, S.},
     TITLE = {On the mathematical foundations of learning},
   JOURNAL = {Bull. Amer. Math. Soc. (N.S.)},
  FJOURNAL = {American Mathematical Society. Bulletin. New Series},
    VOLUME = {39},
      YEAR = {2002},
    NUMBER = {1},
     PAGES = {1--49},
      ISSN = {0273-0979},
   MRCLASS = {68T05 (68T10 91E40 94A20)},
  MRNUMBER = {1864085},
MRREVIEWER = {Andrei Mart\'{\i}nez Finkelshtein},
       DOI = {10.1090/S0273-0979-01-00923-5},
       URL = {https://doi.org/10.1090/S0273-0979-01-00923-5},
}

@article {WSS,
    AUTHOR = {Williamson, R. C. and Smola, A. J. and Sch\"{o}lkopf,
              B.},
     TITLE = {Generalization performance of regularization networks and
              support vector machines via entropy numbers of compact
              operators},
   JOURNAL = {IEEE Trans. Inform. Theory},
  FJOURNAL = {Institute of Electrical and Electronics Engineers.
              Transactions on Information Theory},
    VOLUME = {47},
      YEAR = {2001},
    NUMBER = {6},
     PAGES = {2516--2532},
      ISSN = {0018-9448},
   MRCLASS = {62B10 (47N30 62P15 68T10 94A15 94A17)},
  MRNUMBER = {1873936},
       DOI = {10.1109/18.945262},
       URL = {https://doi.org/10.1109/18.945262},
}

@article {zhou1,
    AUTHOR = {Zhou, D.X.},
     TITLE = {The covering number in learning theory},
   JOURNAL = {J. Complexity},
  FJOURNAL = {Journal of Complexity},
    VOLUME = {18},
      YEAR = {2002},
    NUMBER = {3},
     PAGES = {739--767},
      ISSN = {0885-064X},
   MRCLASS = {68T05 (46E22 62P15)},
  MRNUMBER = {1928805},
MRREVIEWER = {Andrei Mart\'{\i}nez Finkelshtein},
       DOI = {10.1006/jcom.2002.0635},
       URL = {https://doi.org/10.1006/jcom.2002.0635},
}

@article {aronszajn,
    AUTHOR = {Aronszajn, Nachman},
     TITLE = {Theory of reproducing kernels},
   JOURNAL = {Trans. Amer. Math. Soc.},
  FJOURNAL = {Transactions of the American Mathematical Society},
    VOLUME = {68},
      YEAR = {1950},
     PAGES = {337--404},
      ISSN = {0002-9947},
   MRCLASS = {46.0X},
  MRNUMBER = {51437},
MRREVIEWER = {T. H. Hildebrandt},
       DOI = {10.2307/1990404},
       URL = {https://doi.org/10.2307/1990404},
}

@article {Gonzalez2,
    AUTHOR = {Gonzalez, K. and Jord\~ao, T.},
     TITLE = {A close look at the entropy numbers of the unit ball of the
              reproducing {H}ilbert space of isotropic positive definite
              kernels},
   JOURNAL = {J. Math. Anal. Appl.},
  FJOURNAL = {Journal of Mathematical Analysis and Applications},
    VOLUME = {534},
      YEAR = {2024},
    NUMBER = {2},
     PAGES = {17},
      ISSN = {0022-247X,1096-0813},
   MRCLASS = {46E22 (41A46)},
  MRNUMBER = {4693231},
MRREVIEWER = {Haizhang\ Zhang},
       DOI = {10.1016/j.jmaa.2024.128121},
       URL = {https://doi.org/10.1016/j.jmaa.2024.128121},
}

@article {Ingo,
    AUTHOR = {Steinwart, I. and Hush, D. and Scovel, C.},
     TITLE = {An explicit description of the reproducing kernel {H}ilbert
              spaces of {G}aussian {RBF} kernels},
   JOURNAL = {IEEE Trans. Inform. Theory},
  FJOURNAL = {Institute of Electrical and Electronics Engineers.
              Transactions on Information Theory},
    VOLUME = {52},
      YEAR = {2006},
    NUMBER = {10},
     PAGES = {4635--4643},
      ISSN = {0018-9448,1557-9654},
   MRCLASS = {68T05 (46E22 94A12)},
  MRNUMBER = {2300845},
       DOI = {10.1109/TIT.2006.881713},
       URL = {https://doi.org/10.1109/TIT.2006.881713},
}

@book {Michael1,
    AUTHOR = {Ruzhansky, M. and Turunen, V.},
     TITLE = {Pseudo-differential operators and symmetries},
    SERIES = {Pseudo-Differential Operators. Theory and Applications},
    VOLUME = {2},
      NOTE = {Background analysis and advanced topics},
 PUBLISHER = {Birkh\"{a}user Verlag, Basel},
      YEAR = {2010},
     PAGES = {xiv+709},
      ISBN = {978-3-7643-8513-2},
   MRCLASS = {35-02 (35S05 43A77 43A80 43A85 47G30 58J40)},
  MRNUMBER = {2567604},
MRREVIEWER = {Fabio\ Nicola},
       DOI = {10.1007/978-3-7643-8514-9},
       URL = {https://doi.org/10.1007/978-3-7643-8514-9},
}

@book {Michael2,
    AUTHOR = {Fischer, V. and Ruzhansky, M.},
     TITLE = {Quantization on nilpotent {L}ie groups},
    SERIES = {Progress in Mathematics},
    VOLUME = {314},
 PUBLISHER = {Birkh\"{a}user/Springer, [Cham]},
      YEAR = {2016},
     PAGES = {xiii+557},
      ISBN = {978-3-319-29557-2; 978-3-319-29558-9},
   MRCLASS = {22E25 (22E30 35R03 35S05 43A80 46L05)},
  MRNUMBER = {3469687},
MRREVIEWER = {Antoni\ Wawrzy\'{n}czyk},
       DOI = {10.1007/978-3-319-29558-9},
       URL = {https://doi.org/10.1007/978-3-319-29558-9},
}

@article {Delgado,
    AUTHOR = { Delgado, J. and Ruzhansky, M.},
     TITLE = {{$L^p$}-nuclearity, traces, and {G}rothendieck-{L}idskii
              formula on compact {L}ie groups},
   JOURNAL = {J. Math. Pures Appl. (9)},
  FJOURNAL = {Journal de Math\'{e}matiques Pures et Appliqu\'{e}es.
              Neuvi\`eme S\'{e}rie},
    VOLUME = {102},
      YEAR = {2014},
    NUMBER = {1},
     PAGES = {153--172},
      ISSN = {0021-7824,1776-3371},
   MRCLASS = {35S05 (22E30 43A75 47B06)},
  MRNUMBER = {3212252},
MRREVIEWER = {Mohammed\ El A\"{\i}di, Universidad Nacional de Colombia},
       DOI = {10.1016/j.matpur.2013.11.005},
       URL = {https://doi.org/10.1016/j.matpur.2013.11.005},
}

@article {Delgado2,
    AUTHOR = {Delgado, J. and Ruzhansky, M.},
     TITLE = {Schatten classes and traces on compact groups},
   JOURNAL = {Math. Res. Lett.},
  FJOURNAL = {Mathematical Research Letters},
    VOLUME = {24},
      YEAR = {2017},
    NUMBER = {4},
     PAGES = {979--1003},
      ISSN = {1073-2780,1945-001X},
   MRCLASS = {43A75 (22E30 35S05)},
  MRNUMBER = {3723800},
MRREVIEWER = {Sanjiv\ Kumar\ Gupta},
       DOI = {10.4310/MRL.2017.v24.n4.a3},
       URL = {https://doi.org/10.4310/MRL.2017.v24.n4.a3},
}

@book{RKHS,
  title={Reproducing Kernel Hilbert Spaces in Probability and Statistics},
  author={Berlinet, A. and Thomas-Agnan, C.},
  isbn={9781441990969},
  lccn={2003064182},
  url={https://books.google.be/books?id=bX3TBwAAQBAJ},
  year={2011},
  publisher={Springer US}
}

@article {Delgado3,
    AUTHOR = {Delgado, J. and Ruzhansky, M.},
     TITLE = {Schatten--von {N}eumann classes of integral operators},
   JOURNAL = {J. Math. Pures Appl. (9)},
  FJOURNAL = {Journal de Math\'{e}matiques Pures et Appliqu\'{e}es.
              Neuvi\`eme S\'{e}rie},
    VOLUME = {154},
      YEAR = {2021},
     PAGES = {1--29},
      ISSN = {0021-7824,1776-3371},
   MRCLASS = {47G10 (22E30 47B10 58J40)},
  MRNUMBER = {4312282},
MRREVIEWER = {Giuseppe\ Di Fazio},
       DOI = {10.1016/j.matpur.2021.08.006},
       URL = {https://doi.org/10.1016/j.matpur.2021.08.006},
}

@article{Claudinei,
author = {Ferreira, J. and Menegatto, V},
year = {2009},
month = {05},
pages = {61},
title = {Eigenvalues of Integral Operators Defined by Smooth Positive Definite Kernels},
volume = {64},
journal = {Integral Equations and Operator Theory},
doi = {10.1007/s00020-009-1680-3}
}

@book{Cucker, 
place={Cambridge}, 
series={Cambridge Monographs on Applied and Computational Mathematics}, 
title={Learning Theory: An Approximation Theory Viewpoint}, 
publisher={Cambridge University Press}, 
author={Cucker, F. and Zhou, D.X.}, 
year={2007}, 
collection={Cambridge Monographs on Applied and Computational Mathematics}}

@article{Claudinei2,
author = {Ferreira, C.},
year = {2008},
month = { },
pages = { },
title = {Decaimento dos autovalores de operadores integrais gerados por núcleos positivos definidos},
volume = { },
journal = {Dissertação de mestrado, ICMC-USP},
doi = {10.11606/D.55.2008.tde-01042008-091207}
}

@article {Delgado4,
    AUTHOR = {Delgado, J. and Ruzhansky, M.},
     TITLE = {{$L^p$}-bounds for pseudo-differential operators on compact
              {L}ie groups},
   JOURNAL = {J. Inst. Math. Jussieu},
  FJOURNAL = {Journal of the Institute of Mathematics of Jussieu. JIMJ.
              Journal de l'Institut de Math\'ematiques de Jussieu},
    VOLUME = {18},
      YEAR = {2019},
    NUMBER = {3},
     PAGES = {531--559},
      ISSN = {1474-7480,1475-3030},
   MRCLASS = {35S05 (22E30 35R03 47G30)},
  MRNUMBER = {3936641},
MRREVIEWER = {Mattia\ Calzi},
       DOI = {10.1017/s1474748017000123},
       URL = {https://doi.org/10.1017/s1474748017000123},
}

@misc{gonzalez3,
      title={Entropy numbers of Reproducing Hilbert Space of zonal positive definite kernels on compact two-point homogeneous spaces, Arxiv}, 
      journal={J. Fourier Anal. Appl.},
      author={K. Gonzalez and T. Jordão},
      year={2025},
      eprint={2405.08140},
      note={To appear in J. Fourier Anal. Appl.},
      archivePrefix={arXiv},
      primaryClass={math.FA},
      url={https://arxiv.org/abs/2405.08140}, 
}

@article{Moore1,
  author  = {Moore, E. H.},
  title   = {On properly positive Hermitian matrices},
  journal = {Bull. Amer. Math. Soc.},
  volume  = {23},
  pages   = {59},
  year    = {1916}
}

@book{Moorebook,
  author    = {Moore, E. H.},
  title     = {General Analysis},
  series    = {Memoirs of the American Philosophical Society},
  volume    = {1},
  publisher = {American Philosophical Society},
  year      = {1935},
  note      = {Part I, 1935; Part II, 1939}
}

@book{scholkopf_smola_2002,
  author    = {Bernhard Schölkopf and Alexander J. Smola},
  title     = {Learning with Kernels: Support Vector Machines, Regularization, Optimization, and Beyond},
  publisher = {MIT Press},
  year      = {2002},
  address   = {Cambridge, MA},
}

@article{Koeber,
author = {M. Koeber and U. Schäfer},
title = {The unique square root of a positive semidefinite matrix},
journal = {International Journal of Mathematical Education in Science and Technology},
volume = {37},
number = {8},
pages = {990--992},
year = {2006},
publisher = {Taylor \& Francis},
doi = {10.1080/00207390500285867},


URL = { 
    
        https://doi.org/10.1080/00207390500285867
    
    

},
eprint = { 
    
        https://doi.org/10.1080/00207390500285867
    
    

}

}
\end{document}